\documentclass[11pt]{amsart}

\newtheorem{thm}{Theorem}[section]
 \newtheorem{cor}[thm]{Corollary}
 \newtheorem{lem}[thm]{Lemma}
 \newtheorem{prop}[thm]{Proposition}
 \theoremstyle{definition}
 
 \newtheorem{examp}[thm]{Example}
 \theoremstyle{remark}
 
 \theoremstyle{notation}
 
 \numberwithin{equation}{section}

\newcommand{\eps}{\epsilon}

\newcommand{\beq}{\begin{equation}}
\newcommand{\eeq}{\end{equation}}

  \newcommand{\s}{{S}}
 \newcommand{\A}{\mathcal{A}}
 
\newcommand{\f}{\mathcal{F}}

 \newcommand{\bneg}{\mathfrak{b}}

\newcommand{\ngauge}{\mathcal{N}}
\newcommand{\nneg}{\mathfrak{n}}
\newcommand{\gneg}{\mathfrak{m}}
\newcommand{\npos}{\mathfrak{n}^+}
 \newcommand{\g}{\mathfrak{g}}
  \newcommand{\lan}{\mathfrak{l}}

 \newcommand{\Ad}{\textrm{Ad}}
 \newcommand{\ad}{\textrm{ad}}

\newcommand{\gauge}{\mathcal{M}}

 \newcommand{\lop}[1]{\mathfrak{L}(#1)}
\newcommand{\bil}[2]{{\langle #1 | #2\rangle}}

\begin{document}

\title[  $W$-algebras and the equivalence of three reductions ]
 { $W$-algebras and the equivalence of bihamiltonian, Drinfeld-Sokolov and Dirac reductions}

\author[Yassir Dinar]{Yassir Ibrahim Dinar*\\ {\tiny Department of Mathematics and Statistics, Faculty of Science, \\ Sultan Qaboos University, Sultanate of Oman;\\ and The Abdus Salam International Centre for Theoretical Physics (ICTP), Italy.\\ Email: dinar@ictp.it.}}
\thanks{On leave from: Faculty of Mathematical Sciences, University of Khartoum,  Sudan.}




\begin{abstract}
We prove that the classical $W$-algebra associated to a nilpotent orbit in a simple Lie-algebra can be constructed by preforming bihamiltonian,  Drinfeld-Sokolov or Dirac  reductions.
  We conclude that the classical $W$-algebra depends only on the nilpotent orbit but not on the choice of a good grading or an isotropic subspace.
  In addition, using this result we prove again  that the transverse Poisson structure to a nilpotent orbit is  polynomial  and  we
   better clarify  the relation  between classical and  finite $W$-algebras.
\end{abstract}
\maketitle
\tableofcontents
\section{Introduction}

A classical $W$-algebra is a local Poisson bracket on a loop space $\lop M$ of a manifold $M$ where in some  local coordinates
$(u^1,...,u^n)$,  $u^1(x)$ is a Virasoro density
and  $u^i(x),~ i>1$ are primary fields of conformal weights $\eta_i$ \cite{fehercomp}, i.e. they satisfy the identities
\begin{eqnarray}\label{leading terms}
\{u^1(x),u^1(y)\}&=& \eps \delta^{'''}(x-y) +2 u^1(x) \delta'(x-y)+ u^1_x\delta(x-y), \\\nonumber
\{u^1(x),u^i(y)\}&=& (\eta_i+1) u^i(x) \delta'(x-y)+ \eta_i u^i_x \delta(x-y).
\end{eqnarray}
Classical $W$-algebras have a significant role in conformal field theory as their  quantization give  $W$-algebras, i.e.  polynomial  extensions of a particular central extension of the Lie algebra of vector fields on the circle \cite{dickey}.
 They are also associated to integrable hierarchies of partial differential equations of KdV type \cite{gDSh2}. However, we are  interested in classical $W$-algebras  because, possibly after a Poisson  reduction, we can construct  algebraic Frobenius manifolds from
 the dispersionless limit \cite{mypaper},\cite{mypaper1},\cite{mypaper2}.

A wide literature is devoted
to construct examples of classical $W$-algebras within the theory of integrable systems (see \cite{dickey} for some details).
One of the most general and  uniform construction was obtained by Feher et al. in \cite{BalFeh1},
where the authors introduced
a generalization of Drinfeld-Sokolov
reduction that can be performed for any  nilpotent element in simple Lie algebras named
 \textbf{canonical Drinfeld-Sokolov reduction}. This reduction is performed on a standard Lie-Poisson bracket on loop algebra using Dynkin grading, Slodowy slice and
a choice of  a  maximal isotropic subspace (more details are explained below). From the construction it is clear  that nilpotent  elements belonging to the  same nilpotent orbit give equal classical $W$-algebras. By a nilpotent orbit we mean the conjugacy class of nilpotent elements under the adjoint group action. In \cite{BalFeh1}, the authors also argued that canonical Drinfeld-Sokolov reduction is  equivalent  to   Dirac reduction of the Lie-Poisson bracket on Slodowy slice.

Moreover, several attempts have been made to construct  classical $W$-algebras by performing a bihamiltonian reduction
on  Lie-Poisson brackets using the theory of nilpotent orbits.
This was obtained by Casati and Pedroni \cite{CP} to  regular nilpotent orbits in simple Lie algebras via proving the equivalence  between the bihamiltonian and standard Drinfeld-Sokolov reductions. We refer also to  \cite{CASFAL} and \cite{dickey} for  the construction   of classical $W$-algebras associated to  regular nilpotent
orbits in Lie algebras of type $A_n$.

Furthermore, in \cite{mypaper} we obtained  a generalization of the bihamiltonian reduction. This generalization enabled  us to perform bihamiltonian reduction  for any nilpotent orbit in simple Lie algebras. It makes use of the Dynkin grading and  the minimal isotropic subspace.
  In the case of regular nilpotent orbits, our approach made possible  to verify directly that the
  bihamiltonian reduction leads to classical $W$-algebras. For arbitrary nilpotent orbit we proved that the bihamiltonian  reduction is
  equivalent to a generalization of Drinfeld-Sokolov reduction (\cite{mypaper}, section 4).
    Thus to show that the bihamiltonian reduction leads to classical $W$-algebras it is sufficient to prove the equivalence between different types of Drinfeld-Sokolov reductions.

    Actually in this work  we find further results. We prove that the bihamiltonian,  Dirac and Drinfeld-Sokolov reductions  are all equivalent. For a given nilpotent element, we prove that the associated classical $W$-algebra does not depend
    on the choice of a good grading or an isotropic subspace. As a consequence we prove again  that the transverse Poisson structure to a nilpotent orbit is  polynomial  and  we
   better clarify  the relation  between classical and  finite $W$-algebras.

\section{Poisson Geometry and reductions}
In this section we fix some notations  and terminologies. We review our work in \cite{mypaper} and we add some minor results.

A Poisson manifold $M$ is a manifold endowed with a Poisson bracket $\{.,.\}$, i.e. a bilinear skewsymmetric form  on the space of  smooth functions satisfying
 the Leibnitz rule and the Jacobi identity.  Let $M$ be a Poisson manifold with a Poisson bracket $\{.,.\}$.
  Then the   corresponding Poisson tensor $P$ is a linear  map $P: T^*M\to TM$ defined by requiring that
\[ \{F,G\}=\bil {dF}{P ~dG}\]
for any smooth functions $F$ and $G$ on $M$. A smooth function $F$ on $M$ is called a Casimir function,  if it satisfies
\[ P(dF)=\{.,F\}=0.\]

A bihamiltonian manifold $M$ is a manifold endowed with two
Poisson brackets $\{.,.\}_1$ and $\{.,.\}_2$ such that  \[\{.,.\}_\lambda:=\{.,.\}_2+\lambda
\{.,.\}_1\] is a Poisson bracket for any constant $\lambda$. The Jacobi
identity for $\{.,.\}_\lambda$ gives the following equation
\begin{multline}\label{bih:cond}
    \{\{F,G\}_1,H\}_2+\{\{G,H\}_1,F\}_2+\{\{H,F\}_1,G\}_2+\\
    \{\{F,G\}_2,H\}_1+\{\{G,H\}_2,F\}_1+\{\{H,F\}_2,G\}_1=0
\end{multline}
for any smooth functions $F,G$ and $H$ on $M$. It follows from  this equation that the set of all Casimir functions of $\{.,.\}_1$
are closed with respect to $\{.,.\}_2$.

Let  $M$ be a bihamiltonian manifold  with Poisson brackets $\{.,.\}_1$ and $\{.,.\}_2$. Let $P_1$ and $P_2$
denote the corresponding Poisson tensors, respectively.
We assume there is a set
\begin{equation}
\Xi=\{K_1,K_2,...,K_n\}\end{equation} of independent Casimirs of $\{.,.\}_1$  which are closed with respect to $\{.,.\}_2$. For the standard
bihamiltonian reduction \cite{CMP} we assume $\Xi$ to be a complete set of independent Casimirs of $\{.,.\}_1$.
Let us fix a level set $\s$  of $\Xi$ and let $i_s:S \to M$ be the canonical immersion. Then
we consider  the integrable distribution $D$ on $M$ generated by the
Hamiltonian vector fields \begin{equation}
X_{K_i}=P_2(dK_i),\qquad i=1,...,n.
\end{equation}
Let $E$ denotes the distribution induced on $S$ by $D$. We assume the foliation of $E$ on $S$ is regular,
so that $N=S/E$ is a smooth  manifold and $\pi:S\to N$ is a submersion.
 Then applying Marsden-Ratiu reduction theorem \cite{MR}, we get the following result.
 \begin{prop}\cite{mypaper}
The  space $N$ has  a natural bihamiltonian structure $\{.,.\}_2^N$,$ \{.,.\}_1^N$ defined as follows. For any functions $f,g$ on $N$ we have
\begin{eqnarray}
\{f,g\}_2^N\circ \pi&=&\{F,G\}_2\circ i_s\\\nonumber
\{f,g\}_1^N\circ \pi&=&\{F,G\}_1\circ i_s,
\end{eqnarray}
where  $F$ and $G$ are functions on $M$ which extend $f$ and $g$, respectively, and are constant on $D$.
\end{prop}

\subsection{Poisson tensor procedure}

In this section we give a procedure to obtain the reduced bihamiltonian structure, it was introduced for the standard bihamiltonian reduction in \cite{CP}. We assume  that there is a submanifold $Q \subset S$ transverse to $E$, i.e.
\begin{equation}  T_q S=E_q
\oplus T_q Q,\quad\textrm{for all}~q\in Q.
\end{equation}
Then  we have an    isomorphism \[ \Psi:Q \to N\] sending a point to the foliation of $E$ containing that point. The composition  $\Psi^{-1}\circ \pi$ is an inverse of the inclusion map $i_Q: Q \to S$. Hence, the bihamiltonian structure on $N$ can be defined on $Q$  as follows.
For any functions $f,g$ on $Q$ we have
\begin{eqnarray}
\{f,g\}_2^Q&=&\{F,G\}_2\\\nonumber
\{f,g\}_1^Q&=&\{F,G\}_1,
\end{eqnarray}
 where $F, G$ are functions on $M$ extending $f,g$ and constant along $D$. Let  $P_\lambda^Q$ denote the Poisson tensor of $\{.,.\}_\lambda^Q:=\{f,g\}_2^Q+\lambda\{f,g\}_1^Q$.

  \begin{lem} \cite{mypaper}\label{Self:Consistency}
Let  $q\in Q$ and $w \in T_q^*Q$. Then there exists $v\in T_q^* M$
 such that:
\begin{enumerate}
\item $v$ is an extension of $w$, i.e. $(v,\dot{q})=(w,\dot{q})$
for any $\dot{q}\in T_q Q$.
\item $P_\lambda (v) \in T_q Q$.
\end{enumerate} Moreover, the Poisson tensor
$P_\lambda^Q(w)$ is  given by\begin{equation} P_\lambda^Q
w=P_\lambda v \end{equation} for any extension $v$ satisfying
conditions (1) and (2).
\end{lem}
 The previous  lemma leads to a  procedure to calculate the reduced Poisson  bracket. We  refer to  it simply  by \textbf{Poisson tensor procedure}.

\subsection{Bihamiltonian and Dirac reductions}

We show that under further hypothesis the bihamiltonian reduction is equivalent to Dirac reduction.
\begin{cor}  \label{dirac:formula}
In the notations of  lemma \ref{Self:Consistency}, an extension  $v$ of $w$ is unique if and only if
$P_\lambda^Q$ is the Dirac reduction of $P_\lambda$ to $Q$.
\end{cor}
\begin{proof}
We apply Poisson tensor procedure. Let us choose  local coordinates $(q^1,...,q^n)$ on $M$ such that $Q$ is defined by the equations $q^\alpha=0$ for $\alpha=m+1,...,n$. We introduce three types of indices differing by their ranges to simplify  the formulas below; capital letters $I,J,K,...=1,...,n$, small letters $i,j,k,...=1,....,m$ which label the coordinates on the
submanifold $Q$ and Greek letters $\alpha,\beta,\delta,...=m+1,...,n$.
In these notations a covector $w\in T^*Q$ will have the form
\begin{equation}
w=a_i\, dq^i
\end{equation}
and an extension of this covector to $v\in T^*M$ satisfying lemma \ref{Self:Consistency} is given by
\begin{equation}
v=a_I\, d q^I,
\end{equation}
where the coefficients $a_\alpha$'s  are obtained  from requiring that
\begin{equation}\label{eqP}
P_\lambda(v)={P_\lambda}^{IJ} a_J {\partial \over \partial q^I} \in TQ.
\end{equation}
This means that  the  coefficients of ${\partial \over \partial q^\beta}$  equal  $0$ and we get a system of linear equations
\begin{equation}
-{P_\lambda}^{\alpha i} a_i= {P_\lambda}^{\alpha \beta} a_\beta.
\end{equation}
Then   the uniqueness of the extension $v$ is equivalent to the fact that the minor matrix ${P_\lambda}^{\alpha \beta}$ is invertible.
Let ${(P_\lambda)}_{\alpha \beta }$ denote its inverse, then
  \begin{equation}
  a_\beta =- {(P_\lambda)}_{\beta \alpha} {P_\lambda}^{\alpha i} a_i,
  \end{equation}
Substituting this into the formula of $P_\lambda(v)$, we get
 \begin{equation}
P_\lambda(v)=\big( {P_\lambda}^{ij} a_j+{P_\lambda}^{i\beta} a_{\beta}\big){\partial \over \partial q^i}=\big( {P_\lambda}^{ij}-{P_\lambda}^{i\beta} {(P_\lambda)}_{\beta \alpha} {P_\lambda}^{\alpha j}\big)a_j {\partial \over \partial q^i}.
 \end{equation}
Using the identity $P_\lambda^Q(w)=P_\lambda(v)$, we end with  Dirac formula for the reduced Poisson tensor
\begin{equation}\label{dirac formula}
({P_\lambda}^Q)^{ij}={P_\lambda}^{ij}-{P_\lambda}^{i\beta} {(P_\lambda)}_{\beta \alpha} {P_\lambda}^{\alpha j}
\end{equation}
meaning that $P_\lambda^Q$ is the Dirac reduction  of $P_\lambda$ to $Q$.
\end{proof}

We observe that the bihamiltonian reduction   guarantees that
  $P_\lambda^Q$ is linear in $\lambda$ and hence we have a bihamiltonian structure on $Q$.
  This fact is not obvious  when we use  Dirac reduction because  Dirac formula  used to evaluate the reduced Poisson tensor depends on the inverse of a matrix.

\subsection{Local Poisson brackets and Dirac reduction}\label{loc Poiss}

Let  $M$ be a manifold. The loop space $\lop M$ of $M$ is the space of smooth maps from the circle to $M$.
A local Poisson bracket $\{.,.\}$ on $\lop M$ is a Poisson bracket on the space of local  functional on $\lop M$. If we choose local coordinates
$(u^1,...,u^n)$, then $\{.,.\}$  is a finite summation of the form
\begin{eqnarray} \label{genLocPoissBra}\{u^i(x),u^j(y)\}&=&
\sum_{k=-1}^\infty  \{u^i(x),u^j(y)\}^{[k]},\\\nonumber
\{u^i(x),u^j(y)\}^{[k]}&=&\sum_{s=0}^{k+1} A_{k,s}^{i,j}(u(x))
\delta^{(k-s+1)}(x-y),
 \end{eqnarray}
where  $A_{k,s}^{i,j}(u(x))$ are homogenous polynomials in $\partial_x^m
u^r(x)$ of degree $s$ when we  assign
$\partial_x^m u^r(x)$ degree $m$ and $\delta(x-y)$ is the Dirac  delta function defined by
\[\int_{S^1} f(y) \delta(x-y) dy=f(x).\]
In particular, the first terms can be written as follows
\begin{eqnarray}\label{pois:str}
  \{u^i(x),u^j(y)\}^{[-1]} &=& F^{ij}(u(x))\delta(x-y) \\\nonumber
  \{u^i(x),u^j(y)\}^{[0]} &=& F^{ij}_{0}(u(x)) \delta' (x-y)+ \Gamma_k^{ij}(u(x)) u_x^k \delta (x-y)
\end{eqnarray}
where $F^{ij}_0(u)$, $F^{ij}(u)$ and $\Gamma_k^{ij}(u)$ are smooth functions on $M$. It follows from the definition that $F^{ij}(u)$ defines a Poisson bracket on $M$.

Assume we have a local Poisson bracket on the loop space $\lop M$ of a manifold $M$. Let $N\subset M$ be a submanifold of dimension $m$. Then under some assumptions the Poisson bracket can be reduced to $N$ using Dirac reduction. For this end
 we assume $N$ is defined by the equations $u^\alpha=0$ for $\alpha=m+1,...,n$. We introduce three types of indexes; capital letters $I,J,K,...=1,..,n$,
small letters $i,j,k,...=1,....,m$ which parameterize the
submanifold $N$ and Greek letters
$\alpha,\beta,\gamma,\delta,...=m+1,...,n$.

We write the Poisson bracket on $\lop M$ in the form
\[\{u^I(x),u^J(y)\} = \mathbb{F}^{IJ}(u(x))\delta(x-y)\]
where
$\mathbb{F}^{IJ}(u)$ is a matrix differential operator
\begin{equation}
\mathbb{F}^{IJ}(u)=\sum_{k\geq -1}  \sum_{s=0}^{k+1} A_{k,s}^{I,J}(u(x)) {d^{k-s+1} \over dx^{k-s+1}}.
\end{equation}
\begin{prop}
Assume the minor matrix $ \mathbb{F}^{\alpha\beta}(u)$ restricted to $\lop N$ has an inverse $\mathbb{S}_{\alpha\beta}(u)$
 which is a matrix differential  operator of finite order, i.e. a finite sum
\begin{equation}
\mathbb{S}_{\alpha\beta}(u)=\sum_{k= -1}^\infty  \sum_{s=0}^{k+1} B_{k,s}^{\alpha,\beta}(u(x)) {d^{k-s+1} \over dx^{k-s+1}}.
\end{equation}
 Then  Dirac reduction of $\{.,.\}$ to $\lop N$  is well defined and gives a local Poisson structure. The reduced Poisson structure is given by
\[  \{u^i(x),u^j(y)\}^N = \widetilde{\mathbb{F}}^{ij}(u)\delta(x-y)\]
where
\begin{equation}
\widetilde{\mathbb{F}}^{ij}(u)=\mathbb{F}^{ij}(u)-\mathbb{F}^{i\alpha}(u)\mathbb{S}_{\alpha \beta}(u)\mathbb{F}^{\beta j}(u).
\end{equation}
\end{prop}

\begin{proof}
Let $\f$ be a Hamiltonian functional on $\lop {M}$. Then the Hamiltonian flows have the equations
 \begin{equation}\label{DRed1}
 u_t^I=\mathbb{F}^{IJ} \frac{\delta \f}{\delta u^J}
 \end{equation}
   where $\frac{\delta \f}{\delta u^J}$ is the variational derivative of  $\f$ with respect to $u^J(x)$. Following the spirit of \cite{FerNLoc} (see also \cite[chapter 3]{consH}), the Dirac procedure for the reduction of  \eqref{DRed1} to $\lop N$ has the form
 \begin{eqnarray}\label{Dred2}
 u_t^i&=& \mathbb{F}^{iJ} \frac{\delta \f}{\delta u^J}+ \int \{u^i(x),u^\beta(y)\}
 C_\beta(y) dy \\
\nonumber &=& \mathbb{F}^{ij} \frac{\delta \f}{\delta u^j}+
\mathbb{F}^{i\beta}(\frac{\delta \f}{\delta u^\beta}+C_\beta(x)),
\end{eqnarray}
where the Lagrange multiplier  $C_\beta(y)$ is found from the system of linear  equations
\begin{eqnarray}\label{Dred3}
 0=u_t^\alpha &=& \mathbb{F}^{\alpha J} \frac{\delta \f}{\delta u^J}+ \int \{u^\alpha(x),u^\beta(y)\}
 C_\beta(y) dy \\
\nonumber &=& \mathbb{F}^{\alpha j} \frac{\delta \f}{\delta u^j}+
\mathbb{F}^{\alpha\beta}(\frac{\delta \f}{\delta u^\beta}+C_\beta(x)).
\end{eqnarray}
Applying the inverse operator $\mathbb{S}_{\alpha\beta}$, we get
\begin{equation}
\frac{\delta \f}{\delta u^\beta}+C_\beta(x)=- \mathbb{S}_{\beta\alpha}  \mathbb{F}^{\alpha j} \frac{\delta \f}{\delta u^j}.
\end{equation}
Substituting in \eqref{Dred2},
\begin{equation}
 u_t^i=(\mathbb{F}^{ij}- \mathbb{F}^{i\beta}\mathbb{S}_{\beta\alpha}  \mathbb{F}^{\alpha j})\frac{\delta \f}{\delta u^j}.
\end{equation}
Hence, the operator $\widetilde{\mathbb{F}}^{ij}=\mathbb{F}^{ij}- \mathbb{F}^{i\beta}\mathbb{S}_{\beta\alpha}  \mathbb{F}^{\alpha j}$ defines the Poisson bracket of the Dirac reduction of $\{.,.\}$ to $\lop N$.
\end{proof}

We  show the existence of the inverse operator $\mathbb{S}_{\beta\alpha}$ under certain condition.

\begin{prop}\cite{mypaper} \label{dirac-red-for-operators}  In the notations of equation \eqref{pois:str}, if  the minor matrix
$F^{\alpha \beta}$ is nondegenerate on $N$, then the operator $\mathbb{F}^{IJ}$ has an inverse. Moreover, if $F_{\alpha \beta} $ is the inverse matrix of $F^{\alpha \beta}$  and  we write the leading terms of the reduced Poisson bracket on $\lop N$ in the form
\begin{eqnarray}
  \{u^i(x),u^j(y)\}^{[-1]}_N &=& \widetilde{F}^{ij}(u(x))\delta(x-y), \\
  \{u^i(x),u^j(y)\}^{[0]}_N &=& \widetilde{F}^{ij}_0 (u(x))\delta' (x-y)+ \widetilde{\Gamma}_k^{ij}(u(x)) u_x^k \delta (x-y)
\end{eqnarray}
then
\begin{equation}\label{dirac leading term}
\widetilde{F}^{ij}=F^{ij}-F^{i\beta} F_{\beta\alpha} F^{\alpha j},
\end{equation}
\begin{equation}
\widetilde{F}^{ij}_0= F^{ij}_0-F^{i\beta}_0 F_{\beta\alpha}F^{\alpha
j}+F^{i \beta}F_{\beta \alpha} F^{\alpha \varphi}_0 F_{\varphi
\gamma} F^{\gamma j}-F^{i\beta} F_{\beta \alpha} F^{\alpha j}_0
\end{equation}
and
\begin{equation}
\begin{split}
\widetilde{\Gamma}^{ij}_k u_x^k&= \big(\Gamma^{ij}_k - \Gamma^{i\beta}_k  F_{\beta \alpha} F^{\alpha j} + F^{i \lambda} F_{\lambda \alpha} \Gamma^{\alpha \beta}_k F_{\beta \varphi} F^{\varphi j}-F^{i\beta} F_{\beta \alpha} \Gamma^{\alpha j}_k \big) u_x^k\\
&-\big(F^{i\beta}_0 - F^{i\lambda} F_{\lambda \alpha} F^{\alpha \beta}_0 \big)\partial_x(F_{\beta \varphi} F^{\varphi j}).
\end{split}
\end{equation}
The other terms of the reduced Poisson structure  can be found by solving certain recursive equations.
\end{prop}

\begin{cor} \label{Poiss:finite}The Poisson bracket defined on $N$ via
the matrix $\widetilde{F}^{ij}(u)$ equals the Dirac reduction of the Poisson bracket defined on $M$ via the matrix $F^{ij}(u)$.
\end{cor}

\section{Constructing classical $W$-algebra}

We  review some facts about  nilpotent elements in simple Lie algebras.
A good reference for the material in this section  is the book by Collingwood and  McGovern \cite{COLMC}.

We fix  a simple Lie algebra $\g$  over $\mathbb{C}$ and a nilpotent element $f\in \g$. A good grading for $f$ is a $\mathbb{Z}$-grading
\begin{equation}\label{grad1}
\g=\oplus_{i\in \mathbb{Z}} \g_i;~~~~[\g_i,\g_j]\subset \g_{i+j},
\end{equation}
where  \begin{enumerate} \item  $f \in \g_{-2}$, and
\item The map \[\textrm{\ad}~f :\g_j\to \g_{j-2};~~~~~ a \mapsto \ad f(a)=[f,a]\]
is injective for $j\geq 1$ and surjective for $j\leq 1$.
\end{enumerate}
All good gradings for nilpotent elements   are classified
in \cite{ELC}.

 We fix a good grading $\Gamma$ for $f$. Then we choose, by using Jacobson-Morozov theorem,
 a semisimple element  $h$ and a nilpotent element $e \in \g$
such that the set $\A=\{e,h,f\}$ forms an $sl_2$-triple, i.e.
\begin{equation}\label{sl2:relation} [h,e]=2 e,\quad [h,f]=-2
f,\quad [e,f]= h.
\end{equation}
We can assume without loss of generality  that $\A$ is compatible with
$\Gamma$ in the sense that $h\in \g_0$ and $e\in \g_{2}$ \cite{ELC}.

We observe that if $\Gamma'$ denotes the grading on $\g$ defined by
 \[\widehat{\g}_i:=\{x\in \g: \ad~h(x)=i x\},\]
 then it follows from  representation theory of $sl_2$ algebras
that $\Gamma'$ is a good grading for $f$. Such a good grading obtained from a $sl_2$-triple is called  Dynkin grading.
We can map this grading canonically to a weighted Dynkin diagram of $\g$.
It is known that two nilpotent elements are conjugate under the adjoint group action
if and only if they have the same weighted Dynkin diagram. Hence, we  conclude that the construction of classical $W$-algebras, by the methods we will introduce  in next sections, depends only on the nilpotent orbit of $f$.

Let $\bil . . $ denote the Killing form on $\g$. Then there is a natural symplectic bilinear form on $\g_{-1}$ defined by
\begin{equation}
(.,.): \g_{-1} \times \g_{-1} \rightarrow \mathbb{C},~~ (x,y)\mapsto\bil {e}{[x,y]}.
\end{equation}
We use this symplectic structure to fix  an isotropic subspace $\lan\subset \g_{-1}$. Let $\lan'$ denote the symplectic complement of $\lan$ and
 introduce the following  nilpotent subalgebras
\begin{equation}
\gneg:=\lan\oplus \bigoplus_{i\leq -2} \g_i;~~~~~~\nneg:=\lan'\oplus \bigoplus_{i\leq -2} \g_{i}.
\end{equation}
Let $\g_f$ denote  the subspace $\ker (\ad ~f)$ and $\bneg $ denote the orthogonal complement of $\nneg$ under $\bil . .$.
Then from the properties of the good grading we get  \cite{wang}
\begin{equation} \dim \g_f=\dim \g_0+ \dim \g_{-1} ~~\mathrm{ and }~~\g_f\subset \oplus_{i\leq 0} \g_i \subset \bneg.\end{equation}

\begin{lem}\label{slodowy-tansf}
The space $\bneg$ has the following form
\begin{equation}\label{use3}
\bneg=[\gneg,e]\oplus \g_f.
\end{equation}
\end{lem}
\begin{proof}
 We get from  the properties of good grading that  \[0=\bil {[\gneg,\nneg]}{e}=-\bil {\nneg}{[\gneg,e]}\]
which implies that $[\gneg,e]\subset \bneg$. We observe that the properties of $\ad~f$ has its counterpart on $\ad~e$. In particular $\ad ~e:\g_i\to \g_{i+2}$ is injective for $i< 0$. Hence, $\dim[\gneg,e]=\dim \gneg$.
Also, from representation theory of $sl_2$-triples we get \begin{equation}\label{use1} [\gneg,e]\cap \g_f=0.\end{equation}
Computing  the dimension of $\bneg$ we find that \begin{eqnarray}\label{use2}\dim \bneg&=&\dim \oplus_{i\leq 0} \g_i+\dim \g_{1}-\dim \lan'=\dim \oplus_{i\leq 0} \g_i+\dim \lan\\\nonumber &=& \dim \gneg +\dim \g_0+\dim \g_{-1}=\dim [\gneg,e]+ \dim \g_f.\end{eqnarray}
Hence, from \eqref{use1} and \eqref{use2} we get the direct sum \eqref{use3}.
\end{proof}

\subsection{Standard Lie-Poisson structures on loop algebra}
We   define a bihamiltonian structure on the loop algebra $\lop \g$ as follows.
We extend the Killing form on $\g$ to $\lop \g$ by setting
\begin{equation} (u|v)=\int_{S^1}\bil {u(x)}{v(x)} dx,~ u,v \in \lop M.
\end{equation}
We use $(.|.)$ to identify $\lop\g$ with $\lop \g^*$.  We define the gradient $\delta \f (q)$ for a functional $\f$ on $\lop\g$
 to be the unique element in
$\lop\g$ satisfying
\begin{equation}
\frac{d}{d\theta}\f(q+\theta
\dot{s})\mid_{\theta=0}=\int_{S^1}\langle\delta \f|\dot{s}\rangle dx
~~~\textrm{for all } \dot{s}\in \lop\g.
\end{equation}
Then we choose an element  $a\in \g$ which  centralizes the subalgebra $\nneg$, i.e.
\begin{equation}\label{cond on  a}
\nneg\subset\ker\ad \,a.
\end{equation}
Such an element always exists. For example, we can take $a$ to be a homogenous element of the minimal  grading. Finally, we introduce a bihamiltonian structure  $\{.,.\}_2$ and $\{.,.\}_1$ on $\lop \g$, respectively, by means of  Poisson tensors
\begin{eqnarray}\label{bih:stru on g}
P_2(q(x))(v)&=&[ \partial_x+q(x),v(x)]\\\nonumber
P_1(q(x))(v)&=& [a,v(x)],
\end{eqnarray}
for every  $q\in \lop \g$ and $v\in T_q^*\lop \g\cong\lop\g$. It is a well known fact that they define
 a bihamiltonian structure on $\lop \g$ \cite{MRbook}.

 We mention that  $\{.,.\}_2$
  can be interpreted as  the restriction to $\lop \g$ of the Lie-Poisson bracket on the untwisted affine Kac-Moody algebra
 associated to $\g$. In particular, the leading term $\{.,.\}^{[-1]}_2$ defines the Lie-Poisson bracket on $\g$.

\subsection{Generalized Drinfeld-Sokolov reduction}

 We introduce a generalization of Drinfeld-Sokolov reduction by applying Marsden-Weinstein reduction theorem \cite{MR}.

Let us define a gauge action of the adjoint group $\ngauge$ of $\lop {\nneg}$  by
\begin{eqnarray}\label{guage action 1}
q(x)&\rightarrow& \exp( -\ad s(x))[ \partial_x+q(x)]- \partial_x
\end{eqnarray}
where $s(x) \in \lop \nneg$ and $q(x)\in \lop \g$.

\begin{prop}\cite{mypaper}\label{DS as momentum }
The action of $\ngauge$ on $\lop\g$  under the Poisson tensor \[P_\lambda:=P_2+\lambda P_1\] is
 Hamiltonian for all $\lambda$. It admits a momentum map $J$ to be
the projection
\[J:\lop\g\to\lop \npos,\]
 where  $\npos$ is the embedding of $n^*$ in $\g$ under the Killing form. Moreover, $J$ is $\Ad^*$-equivariant.
\end{prop}

We choose $e$ as a regular value of $J$. Since $\bneg$ is the orthogonal complement to $\nneg$, the level set $J^{-1}(e)$ is given by
\beq
S:=\lop\bneg+e.
\eeq

\begin{prop}\label{DS group}
The isotropy group  of $e$ is the adjoint group $\gauge$ of $\lop \gneg$.
\end{prop}
\begin{proof} The isotropy group of $e$ is the subgroup of $\ngauge$ generated by the set
\[G_e=\{s \in \lop \nneg: (\exp (\ad ~s) n, e)=(n,e),~\forall~ n\in \lop \nneg \}.\]
Let $s\in G_e$. Then from the grading properties  we have
\[(\exp (\ad ~s) n, e)=(n,e),~\forall~ n\in \lop \nneg \Leftrightarrow ([s,e],\lop\nneg)=0.\]
The last equality is satisfied if and only if the projection $s_l$ of $s$ to $\lop{\lan'}$ satisfies $([s_l,e],\lop{\lan'})=0$. From the definition this means that $s_l\in \lop \lan$ and therefore $G_e=\lop \gneg$.
\end{proof}

Proposition \ref{DS group} implies, using  Marsden-Weinstein  reduction theorem \cite{MR}, that the space  $S/\gauge$ is a manifold and inherits a Poisson tensor $P_\lambda'$ from  $P_\lambda$.

\subsection{Generalized bihamiltonian reduction}
We perform a bihamiltonian reduction by  considering the set  $\Xi$  of Casimirs of $\{.,.\}_1$ whose gradient  belongs to $\lop \nneg$.
 For example, for any element $b \in \nneg$ we have that
 \[ \f_b(q(x)):= (b|q(x))\] belongs to $\Xi$. Since $\nneg$ is a Lie subalgebra,
 it is easy to verify that $\Xi$ is closed under $\{.,.\}_2$. We take as a level surface the affine subspace
\begin{equation}\s:=\lop \bneg+e
.\end{equation}
Then the  distribution $D$  equals $P_2(\lop \nneg)$. Let $E$ be the restriction of $D$ to $S$, i.e.  $E=P_2(\lop{\nneg})\cap \lop{\bneg}$.

\begin{prop}\label{bih group} The distribution $E$  is given by
\begin{equation} E=P_2(\lop\gneg).
\end{equation}
Moreover, the foliation of $E$ on $S$ is given by the orbits of the adjoint group $\gauge$ of $\lop \gneg$ acting  on $S$ by
\begin{eqnarray}\label{guage action}
s(x)+e&\rightarrow& \exp (-\ad ~m(x))[ \partial_x+s(x)+e]-\partial_x,
\end{eqnarray}
where $m(x) \in \lop \gneg$ and $s(x)\in \lop \bneg$.
\end{prop}
\begin{proof} By definition, $E$ consists of all elements $v\in \lop \nneg$ such that   \begin{equation}
\bil {v_x+[q,v]+[e,v]} {w}=0
\end{equation}
for every  $q\in\lop\bneg$ and $w\in \lop \nneg$. We note that this equation  is satisfied if $v\in \lop {\oplus_{i\leq -2}\g_i}$.
Hence, it is sufficient to assume that $v\in \lop {\lan'}$. But then  $v$ satisfies  the above equation iff  \[\bil {[e,v]}{\lop{\lan'}}=0\]
This implies that $v$ belongs to the symplectic complement $\lop \lan$ of $\lop {\lan'}$. Thus
\begin{equation} E=P_2(\lan)\oplus P_2(\bigoplus_{i\leq -2}\g_i)= P_2(\lop\gneg).
\end{equation}
In proposition \ref{qistranversal}, we prove that the action  \eqref{guage action} is free,  which implies that $E$ is its infinitesimal generator.
\end{proof}

 From this proposition it follows that the space $N=S/\gauge$  is well defined as  it is the orbit space of
 the action \eqref{guage action}.  Hence, we get a bihamiltonian structure $P_1^{N}$ and $P_2^{N}$ on $N$ from $P_1$ and $P_2$, respectively.
 At this point  we already proved the equivalence between Drinfeld-Sokolov and bihamiltonian reductions.
 \begin{thm}\label{thm ds n bih}
  The generalized Drinfeld-Sokolov reduction coincides with the generalized bihamiltonian reduction.
\end{thm}
\begin{proof}
 This  follows directly from propositions \ref{DS group} and \ref{bih group} as in both reductions the reduced  space  is $N=S/\gauge$, where $S=\lop \bneg +e$ and $\gauge$ is  the adjoint group of $\lop \gneg$.
 \end{proof}

  Following the work \cite{DS} and \cite{Pedroni2}, we study the manifold $N$ by introducing  a transverse subspace to the orbits in $S$.
   Slodowy slice is a natural choice of  such transverse subspace since it is coherent with the theory of nilpotent elements.
  It is defined as  the affine  loop subspace
\begin{equation}Q:=e+\lop  {\g_f}\subset S.
\end{equation}

\begin{prop}\label{qistranversal}
The manifold $Q$ is transverse to $E$ on $\s$. Hence, for any element $s(x)+e\in \s$ there is a unique element $m(x) \in
\lop \gneg$ such that
\begin{equation}\label{gauge fix} q(x)+e=\exp (-\ad ~m(x))[ \partial_x+s(x)+e]- \partial_x\end{equation}
belongs to $Q$. The  entries of $q(x)$ give a system of generators for the ring
$R$ of differential polynomials on $S$ invariant under the action \eqref{guage action}.
\end{prop}
\begin{proof} We must prove that for any  $q\in \lop{\g_f}$ and
$\dot{s}\in\lop\bneg$ there are a unique $v\in \lop\gneg$ and a unique $\dot{w} \in
\lop{\g_f}$ such that \begin{equation}\label{rec eq}
\dot{s}=P_2(e+q)(v)+\dot{w}.\end{equation} We write this equation using
the good grading of $\g$. For $t\in\lop \g$, let $t_i$ denote its projection to $\lop {\g_i}$.  Then we can rewrite \eqref{rec eq} as
\begin{equation}
[e,v_{i-2}]+\dot{w_i}=\dot{s_i}-v_i'-\sum_k [q_k,v_{i-k}].
\end{equation}
This gives a linear system of equations which  can be solved recursively because  the map $\ad~ e$ is injective for $i<0$ and we have
\begin{equation}\lop{\g_f}\oplus[e,\lop\gneg]=\lop\bneg\end{equation}
from lemma \ref{slodowy-tansf}. The second part of the proposition can be proved similarly.
\end{proof}

Now we explain  what we call \textbf{Drinfeld-Sokolov method} for calculating the reduced bihamiltonian structure.
We write the coordinates of $Q$ as differential polynomials in the coordinates of $S$ by means of  equation \eqref{gauge fix}
and then apply the Leibnitz rule. If $s^i(x)$ denote the coordinates on $S$, then the Leibnitz rule for $u,v \in R$ have the following form
\begin{equation}\label{caluc DS}
\{u(x),v(y)\}_\lambda={\partial u(x)\over \partial(\partial^m s^i)}\partial_x^m\Big({\partial v(y)\over \partial(\partial^n s^j)} \partial_y^n\big(\{s^i(x),s^j(y)\}_\lambda\big)\Big).
\end{equation}

\subsubsection{Fractional KdV}\label{FKDVex} We  demonstrate  Drinfeld-Sokolov method when  $\g$ is the Lie algebra  $sl_3$  and $f$ is  a minimal nilpotent element. We explain  the different choices of good gradings, isotropic subspaces and first Poisson brackets.  To this end, let us denote $e_{i,j}$ the fundamental $3\times 3$ matrix, i.e. $(e_{i,j})_{s,t}:=\delta_{i,s} \delta_{j,t}$.
We consider the  $sl_2$-triple $\A=\{e,h,f\}$, where $e=e_{1,3}$, $h=e_{1,1}-e_{3,3}$ and $f=e_{3,1}$.
 There are three good gradings compatible with $\A$.
  The following matrices  summarize the degrees   assigned to $e_{i,j}$ by
  these  gradings. The grading $\Gamma_1$ is Dynkin grading.
\begin{equation}
\Gamma_1:=\left( \begin{array}{ccc}
  0 & 1 & 2 \\
  -1 & 0& 1\\
 -2 & -1 &  0
\end{array}\right),~\Gamma_2:=\left( \begin{array}{ccc}
  0 & 0 & 2 \\
  0 & 0& 2\\
 -2 & -2 &  0
\end{array}\right),~\Gamma_3:=\left( \begin{array}{ccc}
  0 & 2& 2 \\
  -2 & 0& 0\\
 -2 & 0 &  0
\end{array}\right)
\end{equation}

 Let us  list some  possible choices for the element $a$ which can be used to define the first Poisson
tensor $P_1$ on $\lop \g$ \eqref{bih:stru on g}. First,  we can take $a=e_{3,1}$
  since it has the minimal degree in all good gradings. We  can also choose $a=e_{3,2}$ (resp. $a=e_{2,1}$)
   since it has the minimal degree in the  grading $\Gamma_2$  (resp. $\Gamma_3$). Moreover,
   we can set $a=e_{2,1}+ e_{3,2}$ (resp. $a=e_{2,1}- e_{3,2}$) when we consider  the grading $\Gamma_1$ and
   fix the isotropic subspace $\lan=\mathbb{C}({e_{2,1}+ e_{3,2}})$  (resp. $\lan=\mathbb{C}({e_{2,1}- e_{3,2}})$).

Under any choice of  a good grading or isotropic subspace, the transverse subspace  $Q$ is the same. We fix for $Q$ the following coordinates. Here we use lower indices for  convenience.
\begin{equation}\label{tFKDV}
q(x)=\left( \begin{array}{ccc}
  q_4(x) & 0 & 1 \\
  q_3(x) & - 2 q_4(x)& 0 \\
  q_1(x) & q_2(x) &  q_4(x)
\end{array}\right).
\end{equation}

Let us consider the grading $\Gamma_1$. We fix the isotropic subspace $\mathbb{C}({e_{2,1}+ e_{3,2}})$ and define $P_1$ by taking  $a= e_{2,1}+ e_{3,2}$. Then  the subspace  $S$ takes the form
\begin{equation}
s(x)=\left( \begin{array}{ccc}
  s_4(x)+s_5(x) & s_6(x) & 1 \\
  s_3(x) & - 2 s_4(x)& -s_6(x)\\
  s_1(x) & s_2(x) &  s_4(x)-s_5(x)
\end{array}\right).
\end{equation}
Equation \eqref{gauge fix} leads to  the following system of generators for the  invariant ring $R$
\begin{eqnarray*}
q_1(x)&=&s_1(x)-\frac{3}{4} s_6^4(x){}+3 s_4(x) s_6^2(x){}-s_2(x) s_6(x)+s_3(x) s_6(x)+s_5^2(x){}-s_5'(x);\\
q_2(x)&=& s_2(x)+s_6(x)^3{}-3 s_4(x) s_6(x)+s_5(x) s_6(x)-s_6'(x);\\
q_3(x)&=& s_3(x)-s_6^3(x){}+3 s_4(x) s_6(x)+s_5(x) s_6(x)-s_6'(x);\\
q_4(x)&=& s_4(x)-\frac{1}{2} s_6^2(x){}.
\end{eqnarray*}
Calculating the reduced Poisson brackets by using the Leibnitz rule \eqref{caluc DS}, the nonzero brackets of $\{.,.\}_1^Q$ are
 \begin{eqnarray}\label{fkdvbih2}
\{q_1(x),q_2(y)\}_1^Q&=&  \frac{3}{2}\,\delta '(x - y)-3\,q_4(x)\delta (x - y);\\\nonumber
\{q_1(x),q_3(y)\}_1^Q&=&   \frac{3}{2}\,\delta '(x - y)+3 q_4(x)\,\delta (x - y);\,  \\\nonumber
\{q_2(x),q_4(y)\}_1^Q&=&  -{1\over 2} \,\delta (x - y);  \\\nonumber
\{q_3(x),q_4(y)\}_1^Q&=& {1\over 2}  \,\delta (x - y),\\\nonumber
\end{eqnarray}
while the nonzero ones of  $\{.,.\}_2^Q$ are
\begin{eqnarray}\label{fkdvbih1}
\{q_1(x),q_1(y)\}_2^Q&=&- \frac{1}{2}\delta ^{'''}(x - y)+ 2\,q_1(x)\,\delta '(x - y) + \partial_x q_1 \delta (x - y); \\\nonumber
\{q_1(x),q_2(y)\}_2^Q&=& \frac{3}{2}\,q_2(x)\,\delta '(x - y) + \frac{1}{2}\,\left( -6\,q_2(x)\,q_4(x) + q_2'(x) \right)\delta (x - y);\\\nonumber
\{q_1(x),q_3(y)\}_2^Q&=&  \frac{3}{2}\,q_3(x)\,\delta '(x - y)+\frac{1}{2}\left( 6\,q_3(x)\,q_4(x) + q_3'(x) \right)\delta (x - y); \\\nonumber
\{q_2(x),q_3(y)\}_2^Q&=&  - \delta ''(x - y)+ \left( q_1(x) - 9q_4(x)^2 - 3q_4'(x)\right)\delta (x - y)\\\nonumber
                      & & -6\,q_4(x)\delta '(x - y) ;\\\nonumber
\{q_2(x),q_4(y)\}_2^Q&=& -\frac{1}{2} q_2(x)\delta (x - y);   \\\nonumber
\{q_3(x),q_4(y)\}_2^Q&=& \frac{1}{2}\,q_3(x) \delta (x - y);\\\nonumber
\{q_4(x),q_4(y)\}_2^Q&=& {1\over 6}\delta' (x - y).
\end{eqnarray}

If we consider the grading $\Gamma_3$ and we define $P_1$ by taking $a=e_{2,1}$, then the space $S$   will take  the  form
 \begin{equation}
s(x)=\left( \begin{array}{ccc}
  s_4(x)+s_5(x) & 0 & 1 \\
  s_3(x) & - 2 s_4(x)& s_6(x)\\
  s_1(x) & s_2(x) &  s_4(x)-s_5(x)
\end{array}\right)
\end{equation}
and  the system of  generators will change to
\begin{eqnarray*}
q_1(x)&=&s_1(x)+s_5(x)^2{}+s_2(x) s_6(x)-s_5'(x);\\
q_2(x)&=& s_2(x);\\
q_3(x)&=& s_3(x)-3 s_4(x) s_6(x)-s_5(x) s_6(x)+s_6'(x);\\
q_4(x)&=& s_4(x).
\end{eqnarray*}
 Calculating  $\{.,.\}_2^Q$  using this  system of generators, we get again the brackets \eqref{fkdvbih1}. This  suggests that the reduced second Poisson bracket is independent of the choice of good grading and isotropic subspace. We prove this result in the next section.

We mention here that the Poisson bracket \eqref{fkdvbih1} is known in the literature as  \textbf{fractional KdV algebra} and  the Poisson bracket \eqref{fkdvbih2} is used in \cite{gDSh2} and \cite{CFMP} to construct an integrable hierarchy.

\subsection{Poisson tensor procedure and Dirac reduction}
Let us apply Poisson tensor procedure to construct $P_\lambda^Q$.

\begin{prop}\label{tensor-proc}
Let $z\in Q$ and  $w\in T_z^*Q$. Then an  extension $v\in T_z^*\lop\g$  of $w$ satisfying the hypothesis of Lemma \ref{Self:Consistency} is unique.
The reduced Poisson tensor in this case  is given by
\begin{equation}
P_\lambda^Q(w)=P_\lambda(v).
\end{equation}
\end{prop}
\begin{proof}
We identify $T_z^*Q\simeq \lop {\g_f}^*$ with $ \lop{\g_e}$ using  the Killing form. Let   $w \in T_z^*Q$. Then a vector $v \in \lop \g$ extends  $w$  if $(w,s)=(v,s)$ for all $s \in \lop {g_f}$. Using the direct sum $\g=[\g,f]\oplus\g_e$, we find that  a vector $v \in \lop \g$ extends  $w$
 if and only if the projection $v_e$ of $v$ to $\lop {\g_e}$ equals $w$. Let us rewrite the condition  $P_\lambda(v) \in T_zQ$  of Lemma \ref{Self:Consistency} under the grading $\Gamma$.
Here for $s\in\lop \g$, we denote $s_i$ its projection to $\lop {\g_i}$. For $i\geq0$, we get  a recursive linear system of equations on the coordinates of $v_{i}$
\begin{equation}\label{recusive}
[v_{i},e]=v_{i+2}'+\lambda [a,v]_{i+2}+\sum_{k\leq 0} [q_k,v_{i+2-k}]
\end{equation}
which can be solved  uniquely since $\ad ~e$ restricted to $\g_i$ is surjective  and the projection of $v_i$ to  kernel $\ad~e$ equals  $(v_e)_i$.
For $i\leq -1$, we have $\g_{i+2}=(\g_f)_{i+2}\oplus[\g_{i},e]$ and we get a recursive linear system of equations on the coordinates of $v_{i}$ by  setting the projection of
\begin{equation}\label{recusive1}
[e,v_{i}]+v_{i+2}'+\lambda [a,v]_{i+2}+\sum_{k\leq 0} [q_k,v_{i+2-k}]
\end{equation}
to $[\g_{i},e]$ equals 0, which  can be solved uniquely as the map $\ad~ e$
restricted to $\g_i$ is injective.
\end{proof}

 Now we are in a position to prove the following theorem.

\begin{thm}\label{bihm-red-is-free}
The reduced second Poisson bracket $\{.,.\}_2^Q$ on $Q$   is independent of the
choice of a good grading  and an isotropic subspace.
\end{thm}
\begin{proof}
We observe that  the calculation of  $P_\lambda^Q$ in  proposition \ref{tensor-proc}  can be done by using any other choice of good grading. This implies that this calculation depends only on the properties of $sl_2$-triples $\{e,h,f\}$. The Poisson bracket  $\{.,.\}_2^Q$ is obtained by setting $\lambda=0$ in the recursive equations \eqref{recusive} and \eqref{recusive1}. This ends the proof.
\end{proof}

We obtain the following theorem by applying  corollary \ref{dirac:formula}.

\begin{thm}\label{dirac-to-bihami}
The Poisson bracket $\{.,.\}_\lambda^Q$ equals  the Dirac reduction of $\{.,.\}_\lambda$ to $Q$. It can be calculated by
using Dirac formulas given in proposition \ref{dirac-red-for-operators}.
\end{thm}

In  \cite{BalFeh1}, the authors proved the following
\begin{thm} When $\Gamma$ is the Dynkin grading  and $\lan$ is a Lagrangian subspace, the Poisson bracket $\{.,.\}_2^Q$ is a classical $W$-algebra.
\end{thm}

Combining this result with theorem \ref{bihm-red-is-free} we get the following

\begin{thm}\label{walg}
 The classical $W$-algebra associated to a nilpotent orbit is independent of the choice of a good grading and an isotropic subspace and it can be calculated equally by using Drinfeld-Sokolov method, Poisson tensor procedure or Dirac formula.
\end{thm}

Let us explain in some details,   how we apply Dirac reduction to find $\{.,.\}_\lambda^Q$. We  fix  a homogenous  basis
$\xi_1,...,\xi_n$ for $\g$ with $\xi_1,...,\xi_m$ a basis for
$\g_f$. Let $\xi^1,...,\xi^n \in \g$ be the dual basis
satisfying \[\langle \xi_i|\xi^j\rangle=\delta_{i}^j.\]
Note  that if $\xi_i\in \g_j$ then $\xi^i\in \g_{-j}$ and $\xi^1,...,\xi^m$ are  a basis for
$\g_e$. We calculate in this basis the   structure constants  and the matrix of the Killing  form
\begin{equation}
[\xi^i,\xi^j]:=c^{ij}_k \xi^k,~~~\bil {[\xi^i,\xi^j]}{a}=c^{ij}_a,~~~g^{ij}=\bil{\xi^i}{\xi^j}.
\end{equation}
Let us consider the following coordinates on  $\lop \g$
\begin{equation}
q^i(z):=\bil{z-e}{\xi^i},~ i=1,\ldots,n.
\end{equation}
Then  matrix differential operator
\begin{equation}\label{bih operator}
\mathbb{F}^{ij}_\lambda=-g^{ij}\partial_x-\sum_{k} c^{ij}_k q^k(x)-\lambda c^{ij}_a
\end{equation}
defines the Poisson brackets
\begin{equation}
\{q^{i}(x),q^j(y)\}_\lambda=\mathbb{F}^{ij}_\lambda\delta(x-y).
\end{equation}
From the construction,  Slodowy slice  $Q$ is defined by $q^\alpha=0$ for $\alpha=m+1,...,n$.  Then we can directly apply  Dirac formulas given in proposition  \ref{dirac-red-for-operators} to find the reduction of $\{.,.\}_\lambda$ to $Q$.

\begin{examp}{(The KdV bihamiltonian structure)}
Let $\g$ be the Lie algebra  $sl_2$ with its standard basis
\begin{equation}
e=\left( \begin{array}{cc}
   0 & 1 \\
   0& 0
\end{array}\right),~h=\left( \begin{array}{cc}
   1 & 0 \\
   0& -1
\end{array}\right), ~f=\left( \begin{array}{cc}
   0 & 0 \\
   1& 0
\end{array}\right).
\end{equation}
For  a point $q\in \lop\g$ we use the notations
\begin{equation}q(x)= q_e(x) e+{1\over 2} q_h(x) h +q_f(x) f\end{equation}
and we define $P_1$ by setting  $a=f$. Then the matrix differential operator on $Q:=e+q_f(x)f$ is given by
\begin{equation}
\mathbb{F}_\lambda^{\alpha,\beta}=\left( \begin{array}{ccc}
   0 & 0& \partial_x \\
   0& 2\partial_x& 0\\
   \partial_x& 0 &0
\end{array}\right)+\left( \begin{array}{ccc}
   0 & 2 (q_f(x)+\lambda)& 0\\
   -2 (q_f(x)+\lambda)& 0& 2\\
   0& -2 &0
\end{array}\right).
\end{equation}
Here, we order the coordinates as $\big(q_f(x), q_h(x),q_e(x)\big)$. The minor matrix operator $\mathbb{F}^{\alpha \beta}_\lambda,~ \alpha,\beta:=2,3$ has the following inverse
\begin{equation}
\mathbb{S}=\left( \begin{array}{cc}
   0 & 0 \\
   0& {1\over 2}\partial_x
\end{array}\right)+ \left( \begin{array}{cc}
   0 & -{1\over 2} \\
   {1\over 2} & 0
\end{array}\right).
\end{equation}
Then apply Dirac formula  to get
\begin{equation}
P_\lambda^Q=-{1\over 2}\partial^3_x+ 2(q_f+\lambda)\partial_x+ q_f
\end{equation}
which gives the bihamiltonian structure associated to  the KdV equation
\begin{equation}
\{q_f(x),q_f(y)\}_2^Q=- \frac{1}{2}\delta ^{'''}(x - y)+ 2(\,q_f(x)+\lambda)\,\delta '(x - y) + \partial_x q_f \delta (x - y).
 \end{equation}
\end{examp}

\section{Conclusions and Remarks}

\subsection{Transverse Poisson structure} Let us consider the leading terms $\{.,.\}_2^{[-1]}$ and ${\{.,.\}}_1^{[-1]}$ of the bihamiltonian structure $\{.,.\}_2$ and $\{.,.\}_1$ on $\lop\g$. In the notations introduced after theorem \ref{walg}, we have
\begin{eqnarray}
\{q^i,q^j\}_2^{[-1]}&=&-\sum_{k}c^{ij}_k q^k,\\\nonumber
{\{q^i,q^j\}}_1^{[-1]}&=&-c^{ij}_a.
\end{eqnarray}
In the same manner as in proposition
\ref{DS as momentum }, we can prove  that the restriction of the   action \eqref{guage action 1} to  the adjoint group of $\nneg$ on $\g$ is
 Hamiltonian  and admits
a momentum map. Taking  $e$ as a regular value, we obtain a bihamiltonian structure  ${\{.,.\}}_1^{Q[-1]}$, ${\{.,.\}}_2^{Q[-1]}$ on Slodowy slice $\widetilde{Q}=e+\g_f$.
 From corollary \ref{Poiss:finite}, this bihamiltonian structure is the leading term of the  bihamiltonian structure $\{.,.\}_\lambda^Q$ on  $Q$.

 The Poisson  structure  ${\{.,.\}}_2^{Q[-1]}$ is known in the literature as the transverse Poisson structure (TPS)
  to the adjoint orbit of $e$. It was originally defined as the Dirac reduction of $\{.,.\}_2^{[-1]}$ to $\widetilde Q$ (see \cite{DamSab} and the references within).  There were many papers  devoted to prove that  the TPS is a polynomial structure. This was not a trivial problem as the method used to calculate the TPS was Dirac formulas and it depends on the inverse of a polynomial matrix. In this paper we proved that, in addition to Dirac formulas, the TPS can be calculated by using Poisson tensor procedure and Drinfeld-Sokolov method. Both lead to a simpler proof  for the polynomiality of the TPS as the former  uses the linear recursive equations obtained in proposition \ref{tensor-proc} and the latter uses the Leibnitz rule \eqref{caluc DS} on differential polynomials.

\subsection{Classical and finite $W$-algebras}

We  mention that  Slodowy slice $\widetilde{Q}$ is associated to the theory of finite $W$-algebras initiated by Premet  \cite{premet}. More precisely,  let $\chi\in \g^*$ be given by
\[ \chi(x)=\bil e x\]
and consider the one dimensional character $\mathbb{C}_{\chi}$ on $\gneg$ given by the restriction of $\chi$.
 Let $U(\g)$ and $U(\gneg)$ be the  universal enveloping algebras of $\g$ and $\gneg$, respectively,
and define the associative algebra   \[Q_\chi:= U(\g)\otimes_{U(\gneg)} \mathbb{C}_{\chi}.\]
Then the finite $W$-algebra is a noncommutative algebra defined as
\beq \label{finit w alg} W_\chi:={\rm End}_{U(\g)}(Q_\chi)^{op}.\eeq
In \cite{GanGin}, Gan and Ginzburg  proved that $W_\chi$  is a quantization of TPS and it is independent of the choice of isotropic subspace, while   Brundan and Goodwin  \cite{bruJon} proved that $W_\chi$ is independent of the choice of  a good grading  (see \cite{wang} and the references within for more details). In this work we  proved a similar argument  for classical $W$-algebras. We hope this will contribute in clarifying more the relation between classical and finite $W$-algebras.

\subsection{ Integrable hierarchies of KdV type}

Let $\{.,.\}_2^Q$ be a classical $W$-algebra associated to a nilpotent element $e$. In this paper we gave
 a procedure to obtain a Poisson bracket $\{.,.\}_1^Q$ such that it forms  with $\{.,.\}_2^Q$ a bihamiltonian structure.
 This Poisson bracket is a reduction of a Poisson bracket defined on $\lop \g$ by means of an
  element $a$ satisfying the following  sufficient condition (see equation \eqref{bih:stru on g}): There exists a good grading $\Gamma$
for $e$ and an isotropic subspace $\lan\subset \g_{-1}$ such that
\begin{equation}\label{con suf1} \nneg:=\lan'\oplus\bigoplus_{i\leq -2} \g_i \subset \ker \ad\,a\end{equation}
where $\lan'$ is the symplectic complement of $\lan$.
 Examples above suggest that this may be a necessary condition as well.  Classifying such elements $a$ may
 help in studying integrable hierarchies associated to  classical $W$-algebras.
  In particular, if $a$  is such that  $a+e$ is regular semisimple then one can obtain an integrable hierarchy
   by using Zakarov-Shabat scheme, i.e. analyzing  the spectrum  of the matrix differential operator
   \[P_\lambda=\partial_x+q(x)+e+\lambda a, ~ q(x)\in \lop \bneg.\]
This  includes  the generalized Drinfeld-Sokolov hierarchy developed in \cite{DS},\cite{gDSh1},\cite{gDSh2} and
\cite{DelFeher}.
 We mention here that in the case of the subregular nilpotent element in the Lie algebra of type $C_3$
  there exist an element $a\in \g$ such that $e+a$ is regular semisimple. Unfortunately,
the sufficient condition \eqref{con suf1} is not satisfied. In other words, the bihamiltonian structure defined
by using this element  $a$ cannot
be reduced to bihamiltonian structure on  Slodowy slice by the methods introduced in this paper.

\subsection{General remark}
It is well known that, under certain assumptions,  from a local bihamiltonian structure on $\lop M$,
 where $M$ is a smooth  manifold, one can construct  a Frobenius structure on  $M$.
 Our main motivation in studying local bihamiltonian structures related to  classical $W$-algebras
 is the classification and construction of algebraic Frobenius manifolds \cite{mypaper},\cite{mypaper1},\cite{mypaper2}. The classification of Frobenius manifolds
  is the first step to classify  local bihamiltonian structures using the concept of central invariants  \cite{DubCentral}.
   In the case of a regular nilpotent element in a simply laced Lie algebra the bihamiltonian  structure obtained from applying standard Drinfeld-Sokolov reduction \cite{DS}
    gives a polynomial Frobenius manifolds and the central invariants are all equal to $1\over 24$.

In a subsequent publication we will consider further examples of Frobenius manifolds and  investigate the central invariants  on bihamiltonian manifolds that are produced by applying the reduction methods introduced in this paper.
\vskip 0.5truecm \noindent{\bf Acknowledgments.}

The author thanks  B. Dubrovin for useful discussions and N. Pagnon for providing the reference \cite{DamSab}. Part of this work was inspired by the "Summer School and Conference in Geometric Representation Theory and Extended Affine Lie Algebras" at the University of Ottawa, organized by the Fields Institute. This work is partially supported by the European Science Foundation Programme ``Methods of Integrable Systems, Geometry,
Applied Mathematics" (MISGAM). The author also like to thank anonymous reviewers who gave corrections and valuable comments that has helped to improve the quality of the manuscript.

\end{document}